\documentclass{amsart}

\usepackage{amssymb,amsmath,amsthm,latexsym}
\usepackage{epsfig}
\usepackage{psfrag}

\newtheorem{Theorem}{\bf Theorem}

\newtheorem{Proposition}[Theorem]{\bf Proposition}
\newtheorem{Corollary}[Theorem]{\bf Corollary}

\newcommand{\be}{\begin{equation}}
\newcommand{\ee}{\end{equation}}

\newcommand{\bt}{\begin{Theorem}}
\newcommand{\et}{\end{Theorem}}
\newcommand{\ei}{\end{itemize}}
\newcommand{\bea}{\begin{eqnarray}}
\newcommand{\eea}{\end{eqnarray}}

\def\C{\mathbb C}

\def\e{{\epsilon}}

\def\k{{\bf k}}

\def\hpic #1 #2 {\mbox{$\begin{array}[c]{l} \epsfig{file=#1,height=#2}\end{array}$}}
\def\wpic #1 #2 {\mbox{$\begin{array}[c]{l} \epsfig{file=#1,width=#2}\end{array}$}}

\begin{document}

\title[\bf Planar algebras and Kuperberg's 3-manifold invariant]
{\bf Planar algebras and Kuperberg's 3-manifold invariant}

\author{Vijay Kodiyalam}
\address{The Institute of Mathematical Sciences, Chennai, India}
\email{vijay@imsc.res.in,sunder@imsc.res.in}
\author{V. S. Sunder}



\begin{abstract} We recapture Kuperberg's numerical
invariant of 3-manifolds associated to a semisimple and cosemisimple
Hopf algebra through a `planar algebra construction'. A  result of
possibly independent interest, used during the proof, which relates duality
in planar graphs and Hopf algebras, is the subject of a final section. 
\end{abstract}

\maketitle

\section{Introduction}
Throughout this paper, the symbol $\k$ will always denote an
algebraically closed field and $H = (H,\mu,\eta,\Delta,\e,S)$ will always denote a semisimple and 
cosemisimple Hopf algebra over $\k$. 
We use $S$ to denote the antipodes of both $H$ and its dual Hopf algebra $H^*$.
The notations $h$ and $\phi$ will be reserved for the unique two-sided
integrals of $H$ and $H^*$ normalised to satisfy $\epsilon(h) =
(dim~H)1_{{\bf k}} = \phi(1_H)$ (in
which case $\phi(h) =  (dim~H)1$).
We will identify $H$ with $H^{**}$ and write the scalar obtained by pairing
$x \in H$ with $\psi \in H^*$ as one of $\psi(x)$, $x(\psi)$, $\langle \psi,x \rangle
$,
or $\langle x,\psi \rangle$. Thus for instance, $\langle \psi,Sx \rangle = \langle x,S\psi \rangle$.

We will need the
formalism of Jones' planar algebras. The basic reference is
\cite{Jns}. A somewhat more leisurely treatment of the basic notions
may also be found in \cite{KdySnd1}. (Mostly, we will follow the latter 
where, for instance, the $*$'s are attached to `distinguished points' on 
boxes rather than to regions.) 

While Vaughan Jones (who introduced planar algebras) mainly looked 
at `$C^*$-planar algebras', which are {\em \`{a} fortiori} defined over $\C$,
we will need to discuss planar algebras over fields possibly different 
from $\C$.
We will, in particular, require some results from \cite{KdySnd2} about 
the planar algebra $P = P(H)$ associated to a semisimple and cosemisimple 
Hopf algebra $H$ over an arbitrary (algebraically closed) field.
(To be entirely precise, we should call it
$P(H,\delta)$, where $\delta$ is a solution in $\k$ of the equation
$\delta^2 = (dim ~H)1$, as we have in \cite{KdySnd2}; but we shall be
sloppy and just write $P(H)$, with the understanding that one choice of
a $\delta$ has been made as above.) In the sequel, we shall freely use
`planar algebra terminology' without any apology; explanations of such
terminology can be found in \cite{KdySnd1} or \cite{KdySnd2}.

This paper is devoted to showing that a `planar
algebra construction', when one works with the planar algebra $P=P(H)$,
yields an alternative construction of Kuperberg's `state-sum
invariant' - see \cite{Kpr} - of a closed 3-manifold associated with $H$.
 
We start with a recapitulation of Kuperberg's
construction, which
involves working with a Heegaard decomposition of the manifold. We
describe Heegaard diagrams in some detail in the short \S 2.
Another short section, \S3, describes our planar algebra
construction.
A long \S4 contains the details of the
verification that the result of our construction agrees with that of
Kuperberg's, and is
consequently an invariant of the manifold.
Given a directed graph $G$ embedded in an oriented 2-sphere,  and a
semisimple cosemisimple Hopf algebra $H$, we associate, in \S5 (which
is self-contained and may be read independently), two elements
$V(G,H)$ and $F(G,H)$ of appropriate tensor powers of $H$. We show
that $V(G,H)$ and $F(G,H^*)$ 
are related via the Fourier transform of the Hopf algebra $H$. 

Our initial verification that Kuperberg's invariant could be obtained by
our planar algebraic prescription depended on the graph-theoretic
result above; what we have presented here is a shorter, cleaner version of
the verification which only uses a special case (Corollary \ref{ngon})
of this result (which latter special case is quite easy to prove
independently).

\section{Kuperberg's invariant of 3-manifolds}

In this section we describe Kuperberg's construction of his invariant.
In addition to Kuperberg's original paper \cite{Kpr}, a very clear
description of the invariant can be found in \cite{BrrWst} which
gives yet another construction.

The only 3-manifolds discussed here will be closed and oriented.
Kuperberg's invariant (which is also defined for 3-manifolds that are
not necessarily closed, though we restrict ourselves to these) is
constructed from a Heegaard diagram of the 3-manifold.
We recall - see \cite{PrsSss} - that a Heegaard diagram consists of
an oriented smooth surface $\Sigma$, say, of genus $g$, and two systems of
smoothly embedded circles on $\Sigma$, which we will denote by $U^1,...,U^g$
and $L_1,...,L_g$ (to conform to Kuperberg's upper and lower circles), such
that each is a  non-intersecting system of curves that
does not disconnect $\Sigma$. (Note that a system of $g$
non-intersecting simple closed curves on a genus $g$ surface will fail
to disconnect it precisely when the complement of the union of
small tubular neighbourhoods of the curves is a 2-sphere with $2g$-holes).
However the
$U$-circles and $L$-circles may well intersect but only transversally.
There is a well-known procedure for constructing a 3-manifold from
such data, and a theorem of Reidemeister and Singer specifies a set of  
moves under which two such Heegaard diagrams determine the same
3-manifold. It is a fact that either (i) reversing the orientation of $\Sigma$
or (ii) interchanging the systems of $U$- and $L$-circles determines the
oppositely oriented 3-manifold.

Consider now a genus $g$ Heegaard diagram
$(\Sigma,U^1,...,U^g,L_1,...,L_g)$.
The computation of Kuperberg's invariant requires a choice of orientation
and base-point
on each of the circles $U^1,...,U^g,L_1,...,L_g$, so fix such a choice.
We assume that none of the base-points is a point of intersection of
a $U$- and an $L$-circle.
Set $K^i_t = U^i \cap L_t, K^i = \coprod_t K^i_t, K_t = \coprod_i
K^i_t, K = \coprod_{i,t} K^i_t$
and let $k^i_t,k^i,k_t,k$ denote their cardinalities
respectively\footnote{$\coprod$ denotes disjoint union.}.
Traverse the circles $L_1$ to $L_g$ in order beginning from their base-points
according to their orientation
and index the points of intersection by the set $I_L = \{(t,p):1 \leq t \leq
g, 1 \leq p \leq k_t\}$, with the lexicographic ordering of $I_L$
agreeing with the order in which the points of $K$ are
encountered.
Refer to this as the `lower numbering' of the points of intersection.
Next, traverse the circles $U^1$ to $U^g$ the same way and
index the points of intersection by the set $I^U = \{(i,j):1 \leq i \leq
g, 1 \leq j \leq k^i\}$, with the lexicographic ordering of $I^U$
agreeing with the order in which the points of $K$ are
encountered.
Refer to this as the `upper numbering' of the points of intersection.
These give bijections $l:I_L \rightarrow K$ and $u:I^U \rightarrow K$.

Consider now the elements
 $\Delta_{k_1}(h) \otimes \cdots \otimes
\Delta_{k_g}(h) \in H^{\otimes k}$
and $\Delta_{k^1}(\phi) \otimes \cdots \otimes
\Delta_{k^g}(\phi) \in (H^*)^{\otimes k}$. 
Also consider, for each $q \in K$, the endomorphism $T_q$ of $H^*$ (or of $H$) defined to be
$id$ or $S$
according as the tangent vectors of the lower and upper
circles at the point $q$, in that order, form a positively or
negatively oriented basis  for the tangent space at $q$ to $\Sigma$.
Kuperberg's invariant is obtained by pairing these off using the bijections
$l$ and $u$ after twisting by the $T_q$.

Here, and elsewhere in this paper,
we will find it convenient to use two bits of Hopf algebra notation:
(i) superscripts
indicate that multiple copies of Haar integrals are being used, while
(ii) subscripts indicate use of our version of the so-called Sweedler
notation for comultiplication - according to which we write, for example,
$\Delta_n(x) = x_1 \otimes \cdots \otimes x_n$
rather than the more familiar
$\Delta_n(x) = \sum_{(x)}x_{(1)} \otimes \cdots \otimes x_{(n)}$
in the interest of notational  convenience.

Thus explicitly, suppose that $c$ and $d$ are the numbers of 
isolated $U$- and $L$-circles repectively in the Heegaard diagram.
Then Kuperberg's invariant is given by the expression:
$$
\delta^{-2g + 2c + 2d} \prod_{q \in K} \langle h^{t(q)}_{p(q)}, T_q \phi^{i(q)}_{j(q)} \rangle
$$
where $t,p$ and $i,j$ are the obvious projection functions on $I_L$ and $I^U$
regarded as functions on $K$ via the $l$ and $u$ identifications respectively.
We may also rewrite this expression as
\begin{equation}\label{later}
\delta^{-2g + 2c} \prod_{t=1}^g  h^{t}\left( \prod_{p=1}^{k_t} T_{l(t,p)} \phi^{i(l(t,p))}_{j(l(t,p))} \right).
\end{equation}
Note that the $\delta^{2d}$ is absorbed into the product as those terms for
which $k_t =0$, each of which gives a $h^t(\epsilon) = \delta^2$.

That this expression is independent of the chosen base-points follows from the traciality of $\phi$ and $h$
on $H$ and $H^*$ respectively while independence of the chosen orientations
follows from the fact of $S$ being an anti-algebra and anti-coalgebra map.
The main result of \cite{Kpr} is that this is a topological invariant
of the 3-manifold determined by the Heegaard diagram
and is, in a sense that is made precise there, complete.
We note that Kuperberg's invariant is a `picture
invariant' in the sense of \cite{DttKdySnd1}.

\section{A planar algebra construction}

In this section, we will describe our method of starting with a
connected, spherical, non-degenerate  planar algebra $P$ with non-zero
modulus $\delta$,  and associating a number to a Heegaard diagram with data
$(\Sigma,U^1,...,U^g,L_1,...,L_g)$ as above.

Associated to such a Heegaard diagram is a certain planar diagram that
conveys the same information.
This is also often called a Heegaard diagram but in order to distinguish
the two, we will refer to the latter picture as a planar Heegaard diagram.
The planar Heegaard digram is obtained from the Heegaard diagram in the
following way. Remove thin tubular neighbourhoods of the $L$-circles from
$\Sigma$ to get
an oriented 2-sphere with $2g$ holes. Now a $U$-circle $U^i$ becomes
either (a) a  simple closed curve on this sphere with holes - in case
$k^i = 0$,  or (b) a collection of $k^i$ arcs with
endpoints on the boundaries of the holes, if $k^i > 0$.

Fix a point on the sphere, and identify its complement with the plane - with
anti-clockwise orientation - and
finally arrive at the associated planar Heegaard diagram, which
consists of the following data:
\begin{enumerate}
\item a set of $2g$ of circles (the boundaries of the tubular
  neighbourhoods of the $L$-circles) that comes in pairs - two circles being
paired off if they come from the same $L$-circle - and denoted
  $L_1^+,L_1^{-},..., L_g^+,L_g^-$ (with $L_i^+$ and $L_i^-$ being
  paired for each $i$, and the choice of which to call $+$ and which
  $-$ being arbitrary);  
\item diffeomorphisms of $L_t^+$ onto $L_t^-$ which
  reverse the orientations inherited by $L_t^\pm$ from the plane;
\item
collections of $k_t$ distinguished points on each of $L_t^+$ and
$L_t^{-}$ - that are points of intersection with the $U$-curves -
which are mapped to one another by the diffeomorphism of (2) above;
\item a collection of curves - which we shall refer to as the strings
  of the diagram -  which are either (a) entire $U$-circles
  which intersect no $L$-circles, or (b) arcs of $U$-curves
  terminating at distinguished points on the $L$-circles.
\end{enumerate}
It is to be noted that the planar Heegaard diagram is specified by the
associated Heegaard diagram together with a `choice of point at infinity'.

From a planar Heegaard diagram we create a planar network in the sense
of Jones.
For this, we will first make a choice of base-points on all the
circles $L^\pm_t$, taking care to ensure that (i) the base-points on
$L^+_t$ and $L^-_t$ correspond under the diffeomorphism (of 2 above)
between $L^\pm_t$, and (ii) the base-points are not on the $U$-curves.

Next, thicken the $U$-curves of the planar Heegaard diagram to black bands.
If the bands are sufficiently thin, no base-point on the $L$-circles will lie
in a black region.
We will refer to the  $L_t^+$ as `positive circles' and the
$L_t^{-}$ as `negative circles'.
Each of the positive and negative circles now has an even number of
distinguished points on its boundary - these being the points of intersection
of the boundaries of the black bands, i.e., the doubled $U$-curves,
with the circles.
For each circle $L_t^\pm$, start from its base-point and move clockwise until
the first band is hit - at a distinguished point - and mark that point
with a $*$. This yields a planar network in Jones'
sense. Call it $N$.

The boxes of this network are the holes bounded by the circles $L^\pm_t$.
There are $2g$ of them with colours $k_1,...,k_g$, each occuring twice,
and we denote these boxes by $B_t^\pm$.
(Recall that $k_t$ is the number of points of intersection of $L_t$
with all the $U$-curves in the original Heegaard diagram.)
Suppose that the boxes of $N$ are ordered as $B_1^+,B_1^-,...,B_g^+,B_g^-$.
The number we wish to associate to the Heegaard diagram is given by
the expression
\be \label{ourkup}
\delta^{-(k_1+k_2+...+k_g)}  Z^P_N(c_{k_1} \otimes
... \otimes c_{k_g})\ee
where $Z^P_N$ is the partition function of the planar network $N$
for the planar algebra $P$ and
$c_{k} \in P_{k} \otimes P_{k}$ is the unique element satisfying
$(id \otimes \tau_{k})((1 \otimes x)c_{k}) = x$ for all $x \in P_{k}$,
and $\tau_{k}$ is the normalised `picture trace' on the $P_{k}$.
The element $c_{k}$ is sometimes referred to as a quasi-basis for
the functional $\tau_{k}$ on $P_{k}$ - see \cite{BhmNllSzl} - and
its existence and uniqueness are guaranteed by the non-degeneracy of
$\tau_k$.
It is true and easy to see that
\begin{equation}\label{ckform}
c_k = \sum_{j \in J} f_j \otimes f^j.
\end{equation}
whenever $\{f_j:j \in J\}$ and $\{f^j:j \in
J\}$ are any pair of bases for $P_k$ which are dual with respect to
the trace $\tau_k$  meaning that
\[\tau_k(f_if^j) = \left\{ \begin{array}{ll} 0 & if ~ i\neq j\\1 & if
    ~i =j \end{array} \right. ~. \]

We will show that when $P =P(H)$, the
expression given by (\ref{ourkup}) agrees with Kuperberg's invariant.

We would like to remark that
the expression given by (\ref{ourkup}) is independent
of the chosen base-points (because $c_{k}$ is
invariant under $Z_R \otimes Z_{R^{-1}}$, where $R$ is the
$k$-rotation tangle) and also independent of the choice of which
circles 
to call positive and which negative, due to the symmetry of $c_{k}$
under the flip (which is an easy consequence of the traciality of $\tau_{k}$).

\section{Concordance with Kuperberg's construction}

Our aim in this section is to show that when $P = P(H)$, the
construction of \S3 yields the same result as that of \S2.

We begin by observing that the construction of the previous section makes
perfectly good sense at the following level of generality.
Let us say that a planar
network is {\em box doubled} if there is given a fixed-point free involution on
the set of its boxes which preserves colours, i.e., its boxes are paired off
with each $k$-box being paired with another such.
Suppose $P$ is a connected, spherical, non-degenerate planar algebra  with
non-zero modulus $\delta$ and $N$ is a box doubled
planar network with $2g$ boxes; let $\sigma \in \Sigma_{2g}$ be any
permutation with the property that the boxes $D_{\sigma(2l-1)}(N)$ and
$D_{\sigma(2l)}(N)$ are paired off, and are of colour $k_l$, say, for
$1 \leq l \leq g$. Then define
\begin{equation}\label{tpn}
\tau^P(N) := \delta^{-(k_1+k_2+...+k_g)}  
Z^P_{\sigma^{-1}(N)}(c_{k_1} \otimes  
\cdots \otimes c_{k_g})
\end{equation}
where, for $\pi \in \Sigma_n$, $\pi(N)$ refers to the network
which is $N$, but with its boxes re-numbered according to $\pi$
- see \cite{KdySnd1}.Thus, again by equation (2.3) of \cite{KdySnd1},
we have
\[\tau^P(N) := \delta^{-(k_1+k_2+...+k_g)}  
Z^P_N(U_\sigma(c_{k_1} \otimes  
\cdots \otimes c_{k_g}) ~,\]
where the notation $U_\sigma$ refers, as in
\cite{KdySnd1}, to the invertible operator  $U_\sigma : \otimes_{i=1}^n
V_i \rightarrow \otimes_{i=1}^n V_{\sigma^{-1}(i)}$ - between $n$-fold tensor
products - defined by
\[ U_\sigma (\otimes_{i=1}^n v_i) = \otimes_{i=1}^n v_{\sigma^{-1}(i)}
~.\]

The motivation for this definition, and in particular for the
normalisation, comes from the $(1+1)$ TQFT of \cite{KdyPtiSnd}.
Symmetry of the $c_{k_j}$ under the flip implies - as in \S 3 - that
the definition $\tau^P(N)$ depends only on $N$, $P$,
and on the pairing between the boxes of $N$, and not on the choice of
the permutation $\sigma$ above.

For the rest of this section, we assume that
\begin{enumerate}
\item $P=P(H)$. (Recall that in this case
$H=P_2$ with non-degenerate trace given by $\tau_2 =
\delta^{-2} \phi$.)
\item $N$ is obtained from a planar Heegaard diagram $D$ - and we
  assume that the choices of $L_t^\pm$ are made in such a way as to
  ensure that the orientation inherited by $L_t^+$ (resp., $L_t^-$)
  from the choice of
  orientation made for $L_t$ in Kuperberg's construction is the
  clockwise (resp., anticlockwise) one.
\item the base points chosen on $L_t^\pm$ to define $N$ correspond to the choices in
  Kuperberg's construction.
\item $N$ has $2g$ boxes $B_1^+,B_1^-,...,B_g^+,B_g^-$ in that order, where the
$B_t^\pm$ have colour $k_t$ and have been paired off as above, with
the boundary of $B_t^\pm$ being identified with $L_t^\pm$. Thus,
  the boxes of $N$   are naturally indexed by
$X = \{(t,\epsilon): 1 \leq t \leq g, \e \in \{+,-\}\}$. (So, we may
choose $\sigma$ to be the identity permutation in the computation of
$\tau^P(N)$.) 
\end{enumerate}

We will proceed to calculate $\tau^P(N)$ in several steps. Our first
step will be to relate $\tau^P(N)$ and $\tau^P(\widetilde{N})$, where
$\widetilde{N}$ is a box doubled planar network that contains only
2-boxes (and is built from $N$).

\begin{figure}[!h]
\begin{center}
\psfrag{k}{\huge $k$}
\psfrag{ck}{\huge $C_k$}
\psfrag{c1}{\huge $C_1$}
\psfrag{c0}{\huge $C_{0_+}$}
\resizebox{9.0cm}{!}{\includegraphics{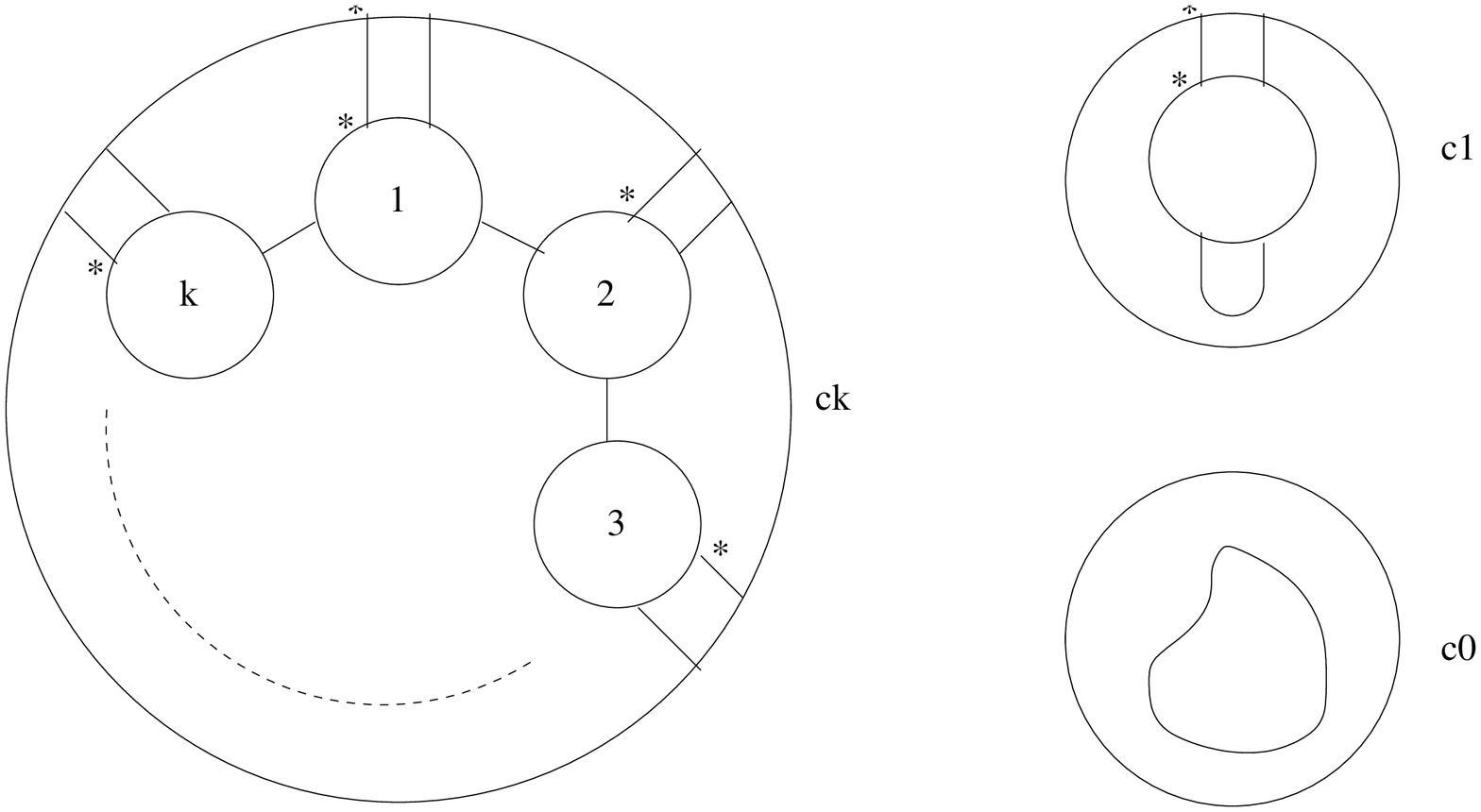}}
\end{center}
\caption{The tangles $C_k$ for $k \geq 2$, $k = 1$ and $k = 0_+$}
\label{fig:tangleck}
\end{figure}

\begin{figure}[!h]
\begin{center}
\psfrag{k}{\huge $k$}
\psfrag{ck}{\huge $C_k^*$}
\psfrag{c1}{\huge $C_1^*$}
\psfrag{c0}{\huge $C_{0_+}^*$}
\resizebox{9.0cm}{!}{\includegraphics{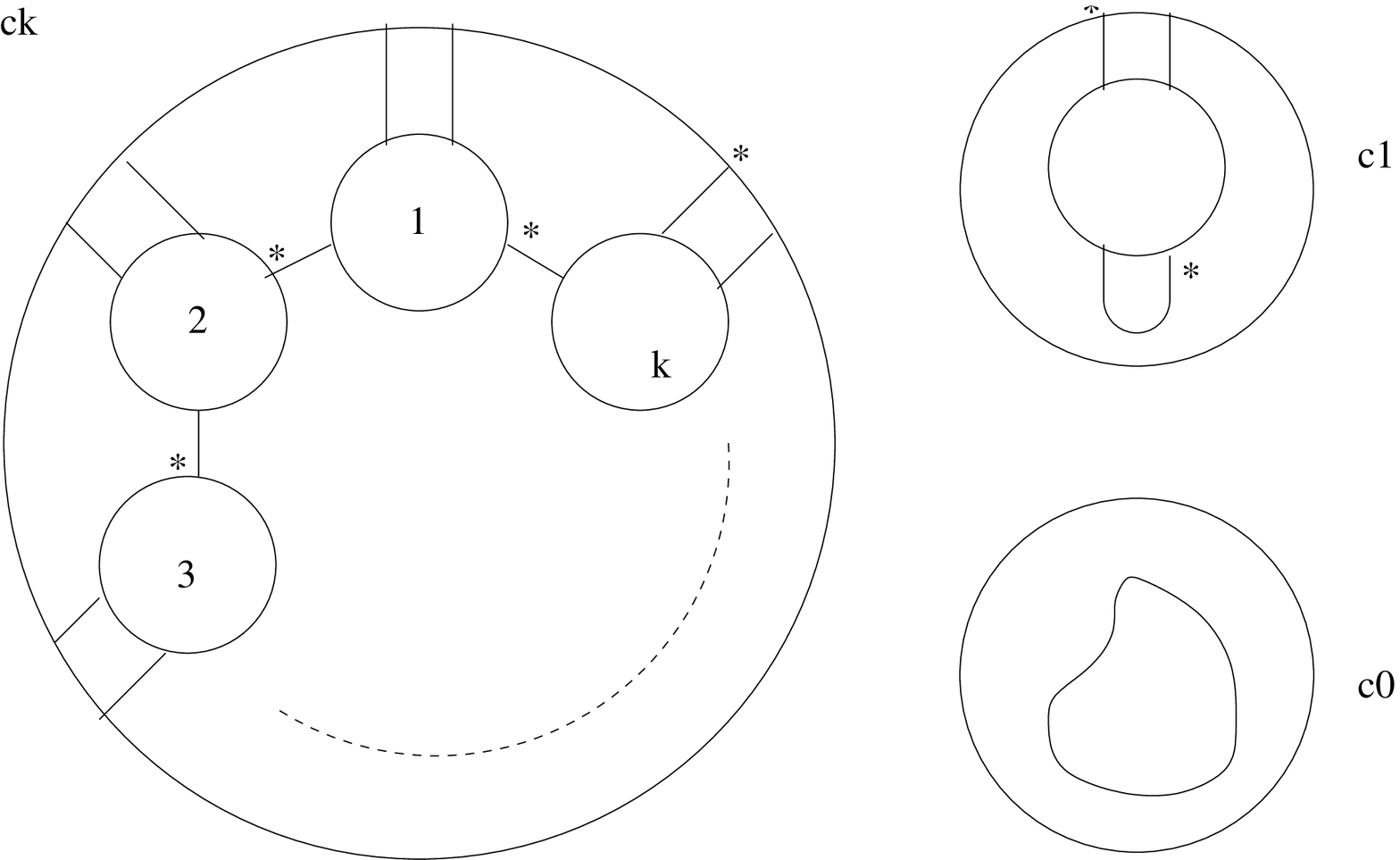}}
\end{center}
\caption{The tangles $C_k^*$ for $k \geq 2$, $k = 1$ and $k = 0_+$}
\label{fig:tangleckst}
\end{figure}

In an obviously suggestive notation, we set $\widetilde{N}$ to be the
planar network defined by 
\[ \widetilde{N} = N \circ_{\{B^\e_t:(t,\epsilon) \in X\}} (\{S(t,\e)\})~,\]
where $S(t,\e)$ is defined to be $C_{k_t}$ or $C_{k_t}^*$ according as
$\e = +$ or $\e =-$, and the tangles $C_k$
are defined in Figure \ref{fig:tangleck} and their adjoint tangles are
illustrated in Figure \ref{fig:tangleckst}.
Note that $\widetilde{N}$ is box doubled, by pairing off the
$p^{th}$ box of $C_{k_t}$ with the $p^{th}$ box of $C_{k_t}^*$.

Our immediate aim is to prove, with the foregoing
notation, that
\begin{equation}\label{tauns}
\delta^{k_1+\cdots+k_g}\tau^P(\widetilde{N}) =
\delta^{2g} \tau^P(N).  
\end{equation}

For this, we begin by noting that in $P_2$, we have
\begin{equation}\label{c2h}
c_2 = h_1 \otimes Sh_2 = Sh_2 \otimes h_1.
\end{equation}

In order to prove equation (\ref{c2h}), note that, for all $x \in H$, we have
\begin{eqnarray*}
(id_H \otimes \frac{1}{n}\phi)((1 \otimes x)(h_1 \otimes Sh_2)) &=&
(id_H \otimes \frac{1}{n}\phi)(h_1 \otimes xSh_2)\\
&=& (id_H \otimes \frac{1}{n}\phi)(h_1x \otimes Sh_2)\\
&=& \frac{1}{n}\phi(Sh_2) h_1x\\
&=& x~;
\end{eqnarray*}
The second identity of equation (\ref{c2h}) is established in similar fashion.

The next step towards proving equation (\ref{tauns}) is to establish
the following identity for $k = 0_+,1,2,\cdots$:
\begin{equation}\label{usigeq}
\delta^2 c_k = (Z_{C_k} \otimes Z_{C_k^*})(U_{\sigma_k}(c_2^{\otimes k}))
\end{equation}
where $\sigma_k \in \Sigma_{2k}$ is the permutation defined by
\[\sigma_k = \left( \begin{array}{ccccccc} 1 & 2 & 3 & 4 & \cdots &
    2k-1 & 2k\\1 & k+1 & 2 & k+2 & \cdots &k & 2k \end{array} \right)~. \]    
 
Note that $U_{\sigma_k}$ maps $H^{\otimes 2k}$ into itself, and we
find from the definition that
\begin{equation}\label{usigform}
U_{\sigma_k}(a(1) \otimes b(1) \otimes a(2) \otimes b(2) \otimes \cdots
\otimes a(k) \otimes b(k) ) = a(1) \otimes \cdots \otimes a(k) \otimes
b(1) \otimes \cdots \otimes b(k) \end{equation}
for any $a(i),b(i) \in H$.

We shall now prove equation (\ref{usigeq}) for $k \geq 2$.
The verification of the equation in the cases $k = 0_+$ and $k=1$  is
easy - and is a consequence of the facts $Z_{C_0}(1) = Z_{C_0^*}(1) =
\delta 1_{0_+}$ and $Z_{C_1} = \epsilon(\cdot)1_1 = (\epsilon \circ S)
(\cdot)1_1 = Z_{C_1^*}$.

We now wish to observe that what was called $X_k$ in Lemma 5 of
\cite{KdySnd2} is nothing but the tangle $C_k \circ_{k} (1^2)$, so that
$Z^P_{X_k}(\otimes_{i=1}^{k-1} a(i)) = Z^P_{C_k}\left((\otimes_{i=1}^{k-1}
a(i)) \otimes 1_H \right)$.
It follows from Lemma 5 of \cite{KdySnd2}, that for $k \geq 1$,  the
LHS of equation (\ref{usigeq}) is given by
\begin{eqnarray}\label{lhs}
\delta^2 c_k &=& \delta^2 \sum_{{\bf i} \in I^{k-1}}  Z_{C_k}(e_{i_1}\otimes
\cdots \otimes  e_{i_{k-1}} \otimes 1) \otimes
Z_{C_k^*}(e^{i_1}\otimes \cdots \otimes   e^{i_{k-1}} \otimes 1) \nonumber\\
&=& (Z_{C_k} \otimes Z_{C_k^*})\left[
  U_{\sigma_k}\left( (\otimes_{j=1}^{k-1}(e_{i_j} \otimes e^{i_j}))
    \otimes (1 \otimes 1) \right)\right] ~~~~\mbox{(by
  eq. (\ref{usigform}))}\nonumber\\  
&=& (Z_{C_k} \otimes Z_{C_k^*})\left[ U_{\sigma_k}\left(c_2^{\otimes
      (k-1)} \otimes (1 \otimes 1)\right) \right] ~~~~~~~~~~\mbox{(by
  eq. (\ref{ckform}))} \nonumber\\
&=& (Z_{C_k} \otimes Z_{C_k^*})\left[ U_{\sigma_k}\left(
    \otimes_{j=1}^{k-1} (h^j_1 \otimes Sh^j_2) \otimes (1 \otimes 1)
  \right) \right] ~~~~\mbox{(by eq. (\ref{c2h}))}\nonumber\\  
&=& \delta^2 Z_{C_k}(h^1_1 \otimes h^2_1 \otimes \cdots \otimes h^{k-1}_1
\otimes 1) \otimes Z_{C_k^*}(Sh^1_2 \otimes Sh^2_2 \otimes \cdots
\otimes Sh^{k-1}_2 \otimes 1) ~.\nonumber\\
&&
\end{eqnarray}

On the other hand, equations (\ref{usigform}) and (\ref{c2h}) imply
that the RHS of equation (\ref{usigeq}) is given by
\[Z_{C_k}(h^1_1 \otimes h^2_1 \otimes \cdots h^k_1) \otimes
Z_{C_k^*}(Sh^1_2 \otimes Sh^2_2 \otimes \cdots Sh^k_2).\]

To proceed further, we need the following consequences of the
so-called `exchange relation' (see \cite{Lnd} and \cite{KdySnd2}) in $P(H)$:
\begin{eqnarray*}
\lefteqn{Z_{C_k}(a(1) \otimes a(2) \otimes \cdots \otimes a(k))}\\
&= &
Z_{C_k}(a(1)Sa(k)_{k-1} \otimes a(2)Sa(k)_{k-2} \otimes
\cdots \otimes a(k-1)Sa(k)_{1} \otimes 1) \\
\lefteqn{Z_{C_k^*}(a(1) \otimes a(2) \otimes \cdots \otimes a(k-1)
  \otimes Sa(k))}\\
 &= &
Z_{C_k^*}(a(k)_{1}a(1) \otimes a(k)_{2}a(2) \otimes
\cdots \otimes a(k)_{k-1}a(k-1)  \otimes 1)
\end{eqnarray*}
for arbitrary $a(1),\cdots,a(k) \in H$. (It is still assumed  that $k$ is
larger than 1.)

We may now deduce that the RHS of equation (\ref{usigeq}) is given by
\begin{eqnarray}\label{rhs}
\lefteqn{Z_{C_k}(h^1_1Sh^k_{k-1} \otimes h^2_1Sh^k_{k-2}\otimes \cdots
h^k_1Sh^k_{1} \otimes 1)}\nonumber\\
&& \otimes Z_{C_k^*}(h^k_kSh^1_2 \otimes
h^k_{k+1}Sh^2_2 \otimes \cdots h^k_{2k-2}Sh^{k-1}_2 \otimes 1) \nonumber\\
&=& Z_{C_k}(h^1_1h^k_k Sh^k_{k-1} \otimes h^2_1h^k_{k+1}Sh^k_{k-2}
\otimes \cdots h^{k-1}_1h^k_{2k-2}Sh^k_{1} \otimes 1) \nonumber\\
&& \otimes
Z_{C_k^*}(Sh^1_2 \otimes
Sh^2_2 \otimes \cdots Sh^{k-1}_2 \otimes 1)\nonumber\\
&&
\end{eqnarray}
where we have used the Hopf algebra fact $xSh_2 \otimes h_1y =
Sh_2 \otimes h_1 xy$ in the last line above. Yet another Hopf algebra
fact guarantees the equality of the right sides of equations
(\ref{lhs}) and (\ref{rhs}); this other (easily established) fact is that
\[h_k Sh_{k-1} \otimes h_{k+1}Sh_{k-2} \otimes \cdots h_{2k-2}Sh_{1}
= \delta^2 ~1^{\otimes (k-1)}~.\]

Now for proving equation (\ref{tauns}), note that
\begin{eqnarray}
\delta^{2g} \tau^P(N) &=& \delta^{2g-(k_1+\cdots +k_g)} Z^P_N(c_{k_1}
\otimes \cdots \otimes c_{k_g}) \nonumber\\
&=& \delta^{-(k_1+\cdots +k_g)}Z^P_N\left[ \otimes_{t=1}^g \left(
  Z^P_{C_{k_t}}\otimes Z^P_{C_{k_t}^*} \right) (U_{\sigma_{k_t}}(c_2^{\otimes
  k_t}))\right] \mbox{(by
  eq. (\ref{usigeq}))} \nonumber\\
&=& \delta^{-(k_1+\cdots +k_g)}Z^P_{\widetilde{N}} (U_{\sigma}(c_2^{\otimes
  k})) \label{ntiln}\\
&=& \delta^{(k_1+\cdots +k_g)}\tau^P(\widetilde{N})~,\nonumber
\end{eqnarray}
where the last step uses the fact that one choice for the permutation
$\sigma \in \Sigma_{2k}$ that is needed in the computation of
$\tau^P(\widetilde{N})$ is given by $\sigma = \coprod_{i=1}^g \sigma_{k_i}$;
and equation (\ref{tauns}) has finally been established.

Next,
note that $K^i$ splits the $U$-circle $U^i$ into $k^i$
strings if $k^i > 0$ or into a single closed string if $k^i=0$. For $(i,j) \in I^U$,
 define $e(i,j)$ to be the string bounded
by $u(i,j-1)$ and $u(i,j)$. (The symbols $l$ and
$u$ refer, of course, to the lower and upper numbering defined in
\S2. Further, we adopt the cyclic convention that $u(i,0)=u(i,k^i)$.)
Orient each string
of the diagram to agree with the choice of orientation of the
$U$-circles in computing Kuperberg's invariant.

We shall use the symbol $E$ to denote the set of non-closed strings of the
diagram $D$ and $C$ to denote the set of closed strings. Thus $|C|$ is
the number of isolated $U$-circles, which was earlier denoted by $c$.
Note that each $e \in E$ comes equipped with the data of
various features of its source and range; specifically, we shall
write:
\begin{itemize}
\item $a(e)$ (resp., $z(e)$) for the point in $K$ at
which the string of the Heegard diagram which corresponds to $e$
  originates (resp., terminates); (these depend only on the original
  Heegaard diagram.)
\item $\alpha(e)$ (resp., $\zeta (e)$) for 1 or 2 according as 
  the string in $D$ which corresponds to $e$ originates (resp.,
  terminates) in a positive or negative box; (these depend on the
  planar Heegaard diagram derived from the original Heegaard diagram.)
\end{itemize}
Note that, by definition,
\be \label{zaeq}
z(e(i,j)) = a(e(i,{j+1})) = u(i,j) ~\forall 1 \leq i \leq g, 1 \leq j \leq k^i~,
\ee
with the convention that $e(i,{k^i+1}) = e(i,1)$. Note also
that the maps
\begin{eqnarray*} \label{zdef}
z,a : E \rightarrow K
\end{eqnarray*}
are bijections and in particular, that $|E| = k$.

We will need to recall the definition and some basic properties of the
Fourier transform map for a semisimple and cosemisimple Hopf algebra.
This is the map $F : H \rightarrow H^*$ defined by $F(x) =
\delta^{-1} \phi_1(x) \phi_2$. The properties that will
be relevant for us are (i) $F \circ F = S $,
(ii) $F \circ S = S \circ F$, (iii)  $F(1) =
\delta^{-1} \phi$ and $F(h) = \delta \epsilon$.
An easily proved Hopf algebra result is:
\begin{equation}\label{ffc2}
 (F \otimes F)(h_1 \otimes Sh_2)
= (\phi_1 \otimes \phi_2).
\end{equation}

We refer the reader to \cite{KdySnd2} for an explanation of the
notations involved and a proof of the following result which
appears as Corollary 10 there.
\begin{Proposition}\label{main}
Let $P = P(H)$ and $Q = P(H^*)$ for a semisimple and cosemisimple
Hopf algebra $H$. Suppose that $N$ is a planar network with $g$ boxes
all of which are 2-boxes. Then:
$$ Z^P_N = Z^Q_{N^-} \circ F^{\otimes g},$$
where both sides are regarded as ${\bf k}$-valued functions on $H^{\otimes g}$.
\end{Proposition}

It follows from Proposition \ref{main}, equation (\ref{ntiln}) and equation (\ref{ffc2}) that
\begin{eqnarray}\label{taupn}
\tau^P(N) &=& \delta^{-(2g + k_1 + \cdots + k_g)} Z^P_{\widetilde{N}}(U_\sigma(c_2^{\otimes k}))\nonumber \\
&=& \delta^{-(2g + k_1 + \cdots + k_g)} Z^Q_{\widetilde{N}^-}
(U_\sigma((\phi_1 \otimes \phi_2)^{\otimes k})).
\end{eqnarray}

We next apply Corollary 3 of \cite{KdySnd2} in order to evaluate 
$Z^Q_{\widetilde{N}^-}
(U_\sigma((\phi_1 \otimes \phi_2)^{\otimes k}))$.
According to this prescription - which was first
outlined in the case of the group planar algebra in \cite{Lnd}, -
given a planar network with only 2-boxes that are labelled by elements of $H$,
its partition function is computed by first replacing
each 2-box labelled by $a$ with
a pair of strands, where the one going through $*$
is labelled $a_1$ and the other $Sa_2$. 
The labels on each loop so formed are read in the order opposite to the
orientation of the loop and $\delta^{-1}\phi$ evaluated on the
product.
The product of these terms over all loops is the required scalar.
We assert that applied to $\widetilde{N}^-$, the number of loops formed is
given by $2g + k + 2c$. 

\begin{figure}[!h]
\begin{center}
\psfrag{x}{1}
\psfrag{y}{2}
\psfrag{1}{(1,1)}
\psfrag{2}{(1,2)}
\psfrag{3}{(1,3)}
\resizebox{6.0cm}{!}{\includegraphics{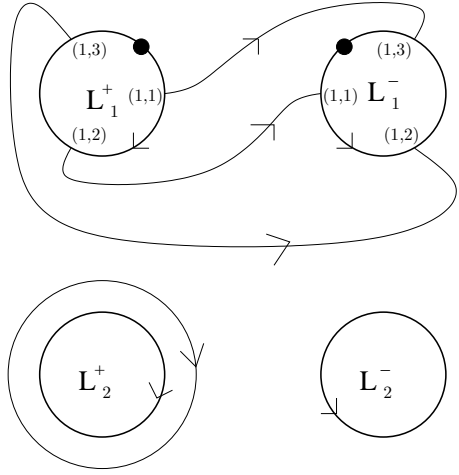}}
\end{center}
\caption{The planar Heegaard diagram for $L(3,1) \# (S^2 \times S^1)$}
\label{fig:connsum}
\end{figure}

For instance consider the planar Heegaard diagram of 
$L(3,1) \# (S^2 \times S^1)$ - the connected sum of the lens space $L(3,1)$
and $S^2 \times S^1$ - shown in Figure \ref{fig:connsum}.
It consists of 2 $U$- and 2 $L$-curves. The $L$ curves have their $\pm$
versions and are shown as dark circles  along with basepoints chosen
on $L_1^\pm$ (the others are irrelevant), while the $U$-curves are shown
by lighter lines. One of the $U$ curves is isolated (the one around $L_2^+$)
while the other breaks up into 3 strings.
The labellings of the points of intersection between the $L$- and $U$-curves
is the `lower numbering'.

\begin{figure}[!h]
\begin{center}
\resizebox{7.0cm}{!}{\includegraphics{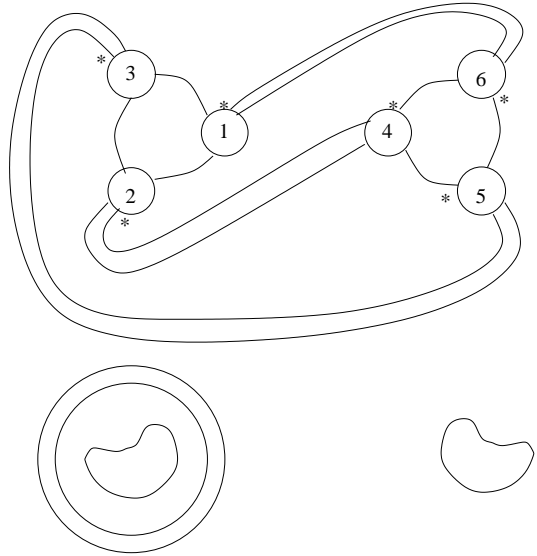}}
\end{center}
\caption{The planar network $\widetilde{N}$ for $L(3,1) \# (S^2 \times S^1)$}
\label{fig:ntilde}
\end{figure}

\begin{figure}[!h]
\begin{center}
\psfrag{p111}{\large $\phi^{(1,1)}_1$}
\psfrag{p121}{\large $\phi^{(1,2)}_1$}
\psfrag{p131}{\large $\phi^{(1,3)}_1$}
\psfrag{p112}{\large $\phi^{(1,1)}_2$}
\psfrag{p122}{\large $\phi^{(1,2)}_2$}
\psfrag{p132}{\large $\phi^{(1,3)}_2$}
\resizebox{8.0cm}{!}{\includegraphics{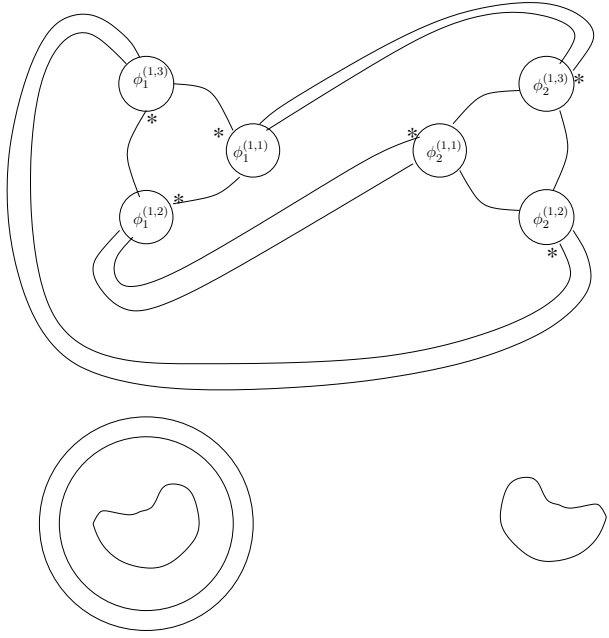}}
\end{center}
\caption{$Z^Q_{\widetilde{N}^-}
(U_\sigma((\phi_1 \otimes \phi_2)^{\otimes 3}))$}
\label{fig:ntildebar}
\end{figure}

\begin{figure}[!h]
\begin{center}
\psfrag{p111}{$\phi^{(1,1)}_1$}
\psfrag{p121}{$\phi^{(1,2)}_1$}
\psfrag{p131}{$\phi^{(1,3)}_1$}
\psfrag{p112}{$S\phi^{(1,1)}_2$}
\psfrag{p122}{$S\phi^{(1,2)}_2$}
\psfrag{p132}{$S\phi^{(1,3)}_2$}
\psfrag{p113}{$\phi^{(1,1)}_3$}
\psfrag{p123}{$\phi^{(1,2)}_3$}
\psfrag{p133}{$\phi^{(1,3)}_3$}
\psfrag{p114}{$S\phi^{(1,1)}_4$}
\psfrag{p124}{$S\phi^{(1,2)}_4$}
\psfrag{p134}{$S\phi^{(1,3)}_4$}
\resizebox{10.0cm}{!}{\includegraphics{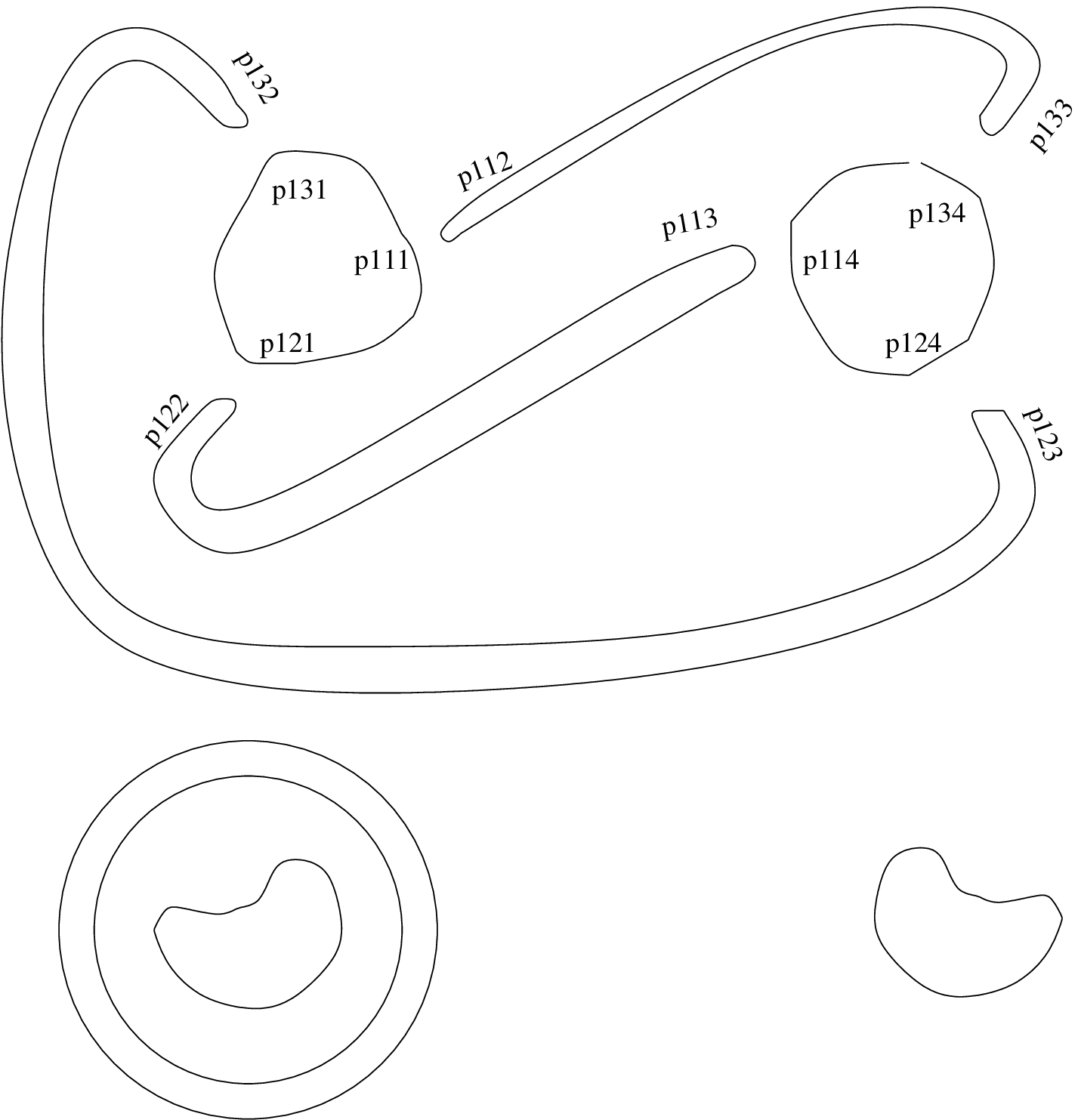}}
\end{center}
\caption{The labelled loops for Figure \ref{fig:ntildebar}}
\label{fig:labloops}
\end{figure}

The planar network $\widetilde{N}$ corresponding to this Heegaard diagram
is shown in Figure \ref{fig:ntilde}.
The planar network $\widetilde{N}^-$ is, by definition, obtained from
$\widetilde{N}$ by moving all the $*$'s anticlockwise by one and
therefore $Z^Q_{\widetilde{N}^-} 
(U_\sigma((\phi_1 \otimes \phi_2)^{\otimes 3}))$ in this example is given
by the labelled planar network in Figure \ref{fig:ntildebar}.
Applying the procedure of Corollary 3 of \cite{KdySnd2} to this labelled planar network
yields the labelled loops as in Figure~\ref{fig:labloops}.
It should now be clear why even in the general case, the number of loops
obtained is $2g + k + 2c$.

Furthermore, a little thought shows that, as in this example,  $Z^Q_{\widetilde{N}^-}
(U_\sigma((\phi_1 \otimes \phi_2)^{\otimes k}))$ is the product of the following
4 types of terms:\\
(a) For each circle of the form $L_t^+$, a term $\delta^{-1}h^{(t,+)}(\prod_{p=1}^{k_t}
\phi^{(t,p)}_1)$,\\
(b) For each circle of the form $L_t^-$, a term $\delta^{-1}h^{(t,-)}(\prod_{p=1}^{k_t}
S\phi^{(t,k_t+1-p)}_4)$,\\
(c) For each closed string in $C$, a multiplicative factor of 
$(\delta^{-1}h(\epsilon))^2 = \delta^2$, and\\
(d) For each non-closed string $e \in E$, a term of the form 
$\delta^{-1}h^e(T_a\phi^{(t_a,p_a)}_{\alpha(e)+1}T_z\phi^{(t_z,p_z)}_{\zeta(e)+1})$, 
where $l^{-1}(a(e)) = (t_a,p_a)$ and $l^{-1}(z(e)) = (t_z,p_z)$ and $T_a$ (resp. $T_z$)
is $S$ or $id$ according as 
$e$ originates (resp. terminates) at a positive or negative box.

Note that (i) since the computation is being done in $Q = P(H^*)$, $h$ and $\phi$
have interchanged roles, as have $1_H$ and $\epsilon$ and (ii) the prescriptions
of (a) and (b) also work for $L_t$'s where $k_t = 0$ with the obvious interpretation
of the empty product.

To summarise, we have seen that
\begin{eqnarray*}
\lefteqn{Z^Q_{\widetilde{N}^-}
(U_\sigma((\phi_1 \otimes \phi_2)^{\otimes k}))}\\
&=&\delta^{-2g+2c-k}
\prod_{t=1}^g h^{(t,+)}(\prod_{p=1}^{k_t}
\phi^{(t,p)}_1)
\prod_{t=1}^g h^{(t,-)}(\prod_{p=1}^{k_t}
S\phi^{(t,k_t+1-p)}_4)
\prod_{e \in E} h^e(T_a\phi^{(t_a,p_a)}_{\alpha(e)+1}T_z\phi^{(t_z,p_z)}_{\zeta(e)+1})\\
&=&\delta^{-2g+2c-k}
\prod_{t=1}^g h^{(t,+)}(\prod_{p=1}^{k_t}
\phi^{(t,p)}_4)
\prod_{t=1}^g h^{(t,-)}(\prod_{p=1}^{k_t}
S\phi^{(t,k_t+1-p)}_3)
\prod_{e \in E} h^e(T_a\phi^{(t_a,p_a)}_{\alpha(e)}T_z\phi^{(t_z,p_z)}_{\zeta(e)})\\
&=&\delta^{-2g+2c-k}
\prod_{t=1}^g h^{(t,+)}(\prod_{p=1}^{k_t}
\phi^{(t,p)}_4)
\prod_{t=1}^g h^{(t,-)}(\prod_{p=1}^{k_t}
S\phi^{(t,k_t+1-p)}_3)
\prod_{e \in E} \phi^{(t_a,p_a)}_{\alpha(e)}(T_ah^e_1) \phi^{(t_z,p_z)}_{\zeta(e)}(T_zh^e_2)
\end{eqnarray*}
where the second equality is a consequence of an application of $\phi_1 \otimes \phi_2
\otimes \phi_3 \otimes \phi_4 = \phi_4 \otimes \phi_1 \otimes \phi_2 \otimes \phi_3$
to each $\phi^{(t,p)}$. We are guilty of a little sloppiness in the
equations above, since actually, $t_a,p_a,t_z,p_z,T_a,T_z$ are all
functions of $e$; 
for instance, $t_a(e) = t(a(e))$ while 
\be \label{TaSTa}
T_a(e) = ST_{a(e)}. 
\ee
(The $T_{a(e))}$ on the right side of the last equation refers to the
$T_q$ used in \S 2.)

Using the relations $Sh = h$ and $h^2 = \delta^2h$, it is easy to see that
\[
\prod_{t=1}^g h^{(t,+)}\left(\prod_{p=1}^{k_t}
\phi^{(t,p)}_4\right)
\prod_{t=1}^g h^{(t,-)}\left(\prod_{p=1}^{k_t}
S\phi^{(t,k_t+1-p)}_3\right)
=
\delta^{2g} \prod_{t=1}^g h^t\left(\prod_{p=1}^{k_t}
\phi^{(t,p)}_3\right),
\]
and therefore we have:
\begin{eqnarray}\label{something}
\lefteqn{Z^Q_{\widetilde{N}^-}
(U_\sigma((\phi_1 \otimes \phi_2)^{\otimes k}))}\\ \nonumber
&=&\delta^{2c-k}
\prod_{t=1}^g h^t\left(\prod_{p=1}^{k_t}
\phi^{(t,p)}_3\right)
\prod_{e \in E} \phi^{(t_a,p_a)}_{\alpha(e)}(T_ah^e_1) \phi^{(t_z,p_z)}_{\zeta(e)}(T_zh^e_2).
\end{eqnarray}

We will next analyse the terms in the product coming from $e \in E$
by grouping together those terms where the $e$'s come from a single $U$-curve.
In other words we write:
\[
\prod_{e \in E} \phi^{(t_a,p_a)}_{\alpha(e)}(T_ah^e_1) \phi^{(t_z,p_z)}_{\zeta(e)}(T_zh^e_2)=
\prod_{\{i: 1 \leq i \leq g, U^i \notin C\}} \prod_{e \subset U^i}\phi^{(t_a,p_a)}_{\alpha(e)}(T_ah^e_1) \phi^{(t_z,p_z)}_{\zeta(e)}(T_zh^e_2)
\]
and for a fixed $i$ such that $U^i \notin C$ (so that $k^i \neq 0$), consider the
product $\prod_{e \subset U^i}\phi^{(t_a,p_a)}_{\alpha(e)}(T_ah^e_1)
\phi^{(t_z,p_z)}_{\zeta(e)}(T_zh^e_2)$. 

Now
$U^i$ comprises of the edges $e(i,j)$ where $1 \leq j \leq k^i$; suppose
$a(e(i,j)) = l(t_{j-1}^i,p_{j-1}^i)$ so that $u(i,j) =z(e(i,j)) = l(t_{j}^i,p_{j}^i)$ (with the
convention that $(t_{0}^i,p_{0}^i) = (t_{k^i}^i,p_{k^i}^i)$).

It follows - from equation (\ref{TaSTa}) - that
\begin{eqnarray*}
\lefteqn{\prod_{e \in U^i}\phi^{(t_a,p_a)}_{\alpha(e)}(T_ah^e_1)
  \phi^{(t_z,p_z)}_{\zeta(e)}(T_zh^e_2)}\\ 
&=&  \prod_{j=1}^{k^i}
  \phi^{(t_{j-1}^i,p_{j-1}^i)}_{\alpha(e(i,j))}(ST_{a(e(i,j))}h^{e(i,j)}_1)
  \phi^{(t_{j}^i,p_{j}^i)}_{\zeta(e(i,j))}(ST_{z(e(i,j))}h^{e(i,j)}_2) ~.
\end{eqnarray*}

After some minor rearrangement, this product may be rewritten as
$$
\prod_{j=1}^{k^i}\phi^{(t_j^i,p_j^i)}_{\zeta(e(i,j))}(ST_{z(e(i,j))}h^{e(i,j)}_2) \phi^{(t_{j}^i,p_{j}^i)}_{\alpha(e(i,j+1))}(ST_{a(e(i,j+1))}h^{e(i,j+1)}_1) 
$$
The definitions show that the $j^{th}$ term of the above product is
$\phi^{(t_j^i,p_j^i)}_1(h_2^{e(i,j)}Sh_1^{e(i,j+1)})$ or $ \phi^{(t_j^i,p_j^i)}_1(h_1^{e(i,j+1)}Sh_2^{e(i,j)})$
according as $e(i,j)$ terminates at a positive or negative
circle. Finally, the product above may be written as:
$$
\prod_{j=1}^{k^i}\phi^{(t_j^i,p_j^i)}_1(ST_{l(t^i_j,p^i_j)}
(h_2^{e(i,j)}Sh_1^{e(i,j+1)})). 
$$

Next, we appeal to Corollary \ref{ngon} $(g - c)$ times - once for
each non-isolated $U^i$. From that corollary, we get:
 $\otimes_{j=1}^{k^i} h_2^{e(i,j)}Sh_1^{e(i,j+1)}
= \delta^{k^i} F^{\otimes k^i} (\Delta_{k^i} \phi^i) = \otimes_{j=1}^{k^i}
\delta F(\phi^i_j)$, which implies that
$$
\prod_{j=1}^{k^i}\phi^{(t_j^i,p_j^i)}_1(ST_{l(t^i_j,p^i_j)}
(h_2^{e(i,j)}Sh_1^{e(i,j+1)})) ~=~
\prod_{j=1}^{k^i}\phi^{(t_j^i,p_j^i)}_1(ST_{l(t^i_j,p^i_j)} 
(\delta F(\phi^i_j))~.
$$
Observe that
$(t,p) = (t^i_j,p^i_j)$ iff $l(t,p) = u(i,j)$ iff $i= i(l(t,p))$ and
$j = j(l(t,p))$.
It now follows from equation (\ref{something}) that
\begin{eqnarray*}
\lefteqn{Z^Q_{\widetilde{N}^-}
(U_\sigma((\phi_1 \otimes \phi_2)^{\otimes k}))}\\
&=&\delta^{2c-k}
\prod_{t=1}^g h^t\left(\prod_{p=1}^{k_t} 
\phi_1^{(t,p)} (S T_{l(t,p)}(\delta F(\phi^{i(l(t,p))}_{j(l(t,p))})))
\phi^{(t,p)}_2\right)\\
&=&\delta^{2c}
\prod_{t=1}^g h^t\left(\prod_{p=1}^{k_t} 
\phi_1^{(t,p)} (S T_{l(t,p)}F(\phi^{i(l(t,p))}_{j(l(t,p))}))
\phi^{(t,p)}_2\right)\\
&=&\delta^{2c}
\prod_{t=1}^g h^t\left(\prod_{p=1}^{k_t} 
\delta FS T_{l(t,p)}F(\phi^{i(l(t,p))}_{j(l(t,p))})
\right)\\
&=&\delta^{2c+k}
\prod_{t=1}^g h^t\left(\prod_{p=1}^{k_t} 
T_{l(t,p)}(\phi^{i(l(t,p))}_{j(l(t,p))})
\right)
\end{eqnarray*}

Finally, a perusal of equations (\ref{later}) and (\ref{taupn}) completes
the verification that Kuperberg's invariant is indeed given by $\tau^P(N)$.

\section{On spherical graphs and Hopf algebras}

Throughout this section, the symbol $G$ will denote an oriented graph
embedded on an oriented smooth sphere $S^2$.
Thus $G$ comprises of a finite subset $V \subset S^2$ of vertices and
a finite set $E$ 
of edges. We regard an edge $e \in E$ as a smooth map from the unit
interval $I$ to $S^2$ such that $e(0),e(1) \in V$ and such that $e$ is
injective except possibly that $e(0) = e(1)$.
Two (images of) distinct edges do not intersect except possibly at vertices.
Thus multiple edges and self-loops are allowed.
An edge $e$ is regarded as being oriented from $e(0)$ to $e(1)$.
We regard $G$ as the subset of $S^2$ given by the union of its edges and
isolated vertices, if any.
By a face of $G$, we mean a connected component of the complement of $G$
in $S^2$.

We will use the terms anticlockwise and clockwise to stand for ``agreeing
with the orientation of" and ``opposite to the orientation of" $S^2$ 
respectively.
If $u$ is the direction of the oriented edge $e$ at a point $p$, and
if $v$ is a perpendicular direction such that $\{u,v\}$ is positively
(resp., negatively) oriented (according to the orientation of the
underlying $S^2$), we shall call the points near $p$ on the side
indicated by $v$ as the 'left' (resp., `right') of the edge $e$.

We digress now with a discussion of tensor products of indexed families
of vector spaces. We consider only finite indexing sets.
For a family $\{V_q :q \in K\}$  of vector spaces (over some
field $\k$), which is indexed by the
finite set $K$, we define $\otimes_{q \in K} V_q$ to be the quotient
of the vector space, with basis consisting of functions $f : K \rightarrow
\coprod_{q \in K} V_q$ such that $f(q) \in V_q$ for all $q \in
K$, by the subspace
spanned by
\[\{f - \alpha_1f_1 - \alpha_2f_2: \exists q_0 \in K \mbox{
  such that}  f(q)= \left\{ \begin{array}{ll} f_1(q)=f_2(q) & if ~q \neq q_0\\
\alpha_1f_1(q)+ \alpha_2f_2(q) & if ~q=q_0 \end{array} \right.\}~.\]

We denote the image in $\otimes_{q \in K} V_q$ of the function $f$
by $\otimes_{q \in K} f(q)$.
If $\{T_q : V_q \rightarrow W_q\}_q$ is an indexed family of vector space maps,
there is a natural induced map $\otimes_{q \in K} T_q: \otimes_{q
  \in K} V_q \rightarrow \otimes_{q \in K} W_q$.

In the important special case of this indexed tensor product when $V_q
= V$ for all $q \in K$, we will also  
denote $\otimes_{q \in K} V_q$ by $V^{\otimes K}$.
We adopt a similar convention for tensor product of vector space maps.

Note that if $K = \{1,2,...,k\}$, then $\otimes_{q \in K} V_q$ can be
naturally identified with $\otimes_{q=1}^k V_q = V_1 \otimes \cdots
\otimes V_k$, and in particular, we will write $V^{\otimes K} = V^{\otimes k}$.
More generally, if $K$ is a totally ordered finite set with $|K| = k$,
then $V^{\otimes K}$ can be naturally identified with $V^{\otimes k}$.
Even more generally, a bijection,
say $\theta$,
from a set $L$ to a set $K$, induces a functorial isomorphism, which we will
denote by $\tilde{\theta}$, from $\otimes_{l \in L}V_{\theta(l)}$ to
$\otimes_{k \in K}V_k$, and in particular from
$V^{\otimes L}$ to $V^{\otimes K}$.
In the sequel, we will use without explicit mention, the canonical
identifications
\begin{eqnarray*}
V^{\otimes (\coprod_{i \in I} K_i)} &\sim& \otimes_{i \in I} V^{\otimes K_i}\\
(V^{\otimes K})^{\otimes L} &\sim& V^{\otimes(L \times K)}.
\end{eqnarray*}

To the pair $(G,H)$ (of a graph and a Hopf algebra), we shall associate
two elements of $H^{\otimes E}$. One of these is computed using
the faces of $G$ and is denoted by $F(G,H)$ and the other is
computed using the vertices of $G$ and is denoted by $V(G,H)$.
The main result of the section relates $F(G,H^*)$ and $V(G,H)$.

\begin{figure}[!h]
\begin{center}
\resizebox{10.0cm}{!}{\includegraphics{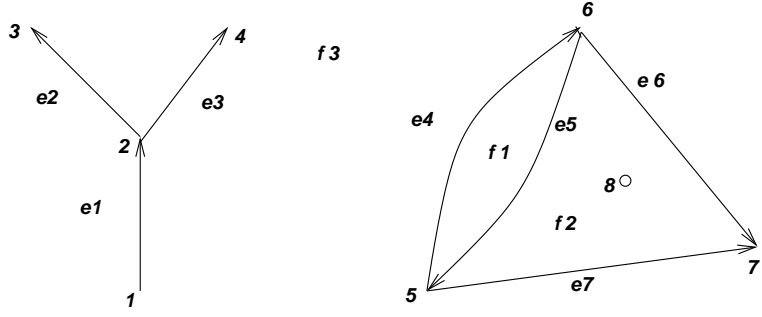}}
\end{center}
\caption{The graph $G$}
\label{fig:eg}
\end{figure}

We will make use of the example illustrated in Figure \ref{fig:eg} -
of a directed graph $G$ 
with eight vertices and three faces - with multiple edges ($e4$  and
$e5$ between vertices 5 and 6) and an isolated vertex (vertex 8) - to
clarify our definitions.

Let $D(V)$ denote the set $E \times \{0,1\}$.
For a vertex $v \in V$, let $D_v$ denote the set $\{(e,i) \in D(V) : e(i) = 
v\}$ and let $d_v$ denote its cardinality which is the degree of $v$.
Consider an enumeration of $D_v$  in 
clockwise order around the vertex $v$. 
This is, of course, determined once one of the edges at $v$ is chosen as the
first.\footnote{For our example, the sets $D_v$, with their
 elements listed in a possible order, are: $D_1 = \{(e1,0)\},
D_2 = \{(e1,1),(e2,0),(e3,0)\},
D_3 = \{(e2,1) \},
D_4 = \{(e3,1) \},
D_5 = \{(e4,0),(e5,1),(e7,0) \},
D_6 = \{(e4,1),(e6,0),(e5,0) \},
D_7 = \{(e6,1),(e7,1) \},
D_8 = \emptyset $.}
Denote this bijection by $\theta_v : \{1,\cdots,d_v\} \rightarrow D_v$.
Note that $D(V)$ is the disjoint union of $D_v$ as $v$ varies over $V$ and 
consider $\otimes _{v \in V}
\widetilde{\theta_v}(\delta^{-1}\Delta_{d_v}(h)) \in H^{\otimes 
D(V)}$.
The traciality of $h$ implies that this element is independent of the
choice of clockwise ordering of the edges around each vertex.

Now consider the map $\mu \circ (id \otimes S): H^{\otimes \{0,1\}} = 
H^{\otimes 2} \rightarrow H$
and the tensor product map $\otimes_{e \in E} (\mu \circ (id \otimes S)) :
H^{\otimes D(V)} = H^{E \times \{0,1\}} \rightarrow H^E.$ 
Define $V(G,H)$ to be the image under this map of $\otimes _{v \in V} 
\widetilde{\theta_v}(\delta^{-1}\Delta_{d_v}(h))$. Explicitly, we have
\be \label{vghdef}
V(G,H) = \delta^{\rho(G)}\otimes_{e \in E} h^{s(e)}_{m(e)}
Sh^{r(e)}_{n(e)} 
\ee
where (i) $\rho(G) =  - |V| + 2|\{v \in V: d_v = 0\}|$;\footnote{The
    reason for the correction term `$+2|\{v:d_v=0\}|$' is that
    $\Delta_0(h) = \e(h) = n = \delta^2$.} and
(ii)  $s,r,m,n$ are functions defined on $E$ and with appropriate
ranges, so that $(e,0)$ is the $m(e)$-th
element of $D_{s(e)}$ while $(e, 
1)$ the $n(e)$-th element of  of $D_{r(e)}$, for any edge $e \in E$.
(Thus, for example, $s,r :E \rightarrow V$ are the `source' and
`range' maps.)

For our example, $V(G,H) \in
H^{\otimes 7}$ - since there are 7 edges; the prescription unravels to yield
\be \label{vgheg}
V(G,H) = \delta^{-6} \left(h^1Sh^2_1 \otimes h^2_2 Sh^3 \otimes
  h^2_3Sh^4 \otimes h^5_1Sh^6_1 \otimes h^6_3Sh^5_2 \otimes
  h^6_2Sh^7_1 \otimes h^5_3Sh^7_2\right)~.\ee

A similar construction using the faces yields $F(G,H)$. 
For this, begin with the set $D(F) = E \times \{l,r\}$.
Consider a pair $(f,c)$ where $f$ is a face of $G$ and $c$ is a component
of the boundary of $f$.
By $\widetilde{F}$, we will refer to the set of all such pairs.
(This set is the set `dual' to the vertex set $V$ in case the graph
$G$ is disconnected.)
Let $D_{(f,c)} = \{(e,d) \in D(F) : e(t) \in c ~{\rm ~for~all~}~ t \in [0,1]$ 
~ and there exist points in $f$ sufficiently close to $c$ where the
orientation agrees or disagrees with the orientation of $e$ according
as $d$ is $l$ or $r$ $\}$. 
We pause to explain this mouthful of a definition. A pair consisting
of an edge $e$ and a direction $d$ is put into $D_{(f,c)}$ exactly
when the image of the edge is part of $c$ and some parts of $f$ lie
to the left or right of $e$ according as $d$ is $l$ or $r$.
Note that it is quite possible for points of $f$ to lie on
both sides of the image of $e$.
Set $d_{(f,c)}$ to be the cardinality of $D_{(f,c)}$.

In our example, there are three faces $f1,f2,f3$, and these
boundaries have  $1, 2 \mbox{ and } 2$ components respectively, and we
have 
\[\widetilde{F} = \{\tilde{f}1 =(f1,4\bar{5}), \tilde{f}2 =
(f2,\bar{5}6\bar{7}), \tilde{f}3 = (f2,\cdot), \tilde{f}4 =
(f3,12\bar{2}3\bar{3}\bar{1}), \tilde{f}5 = (f3,\bar{4}7\bar{6})\}
~,\] 
with the  notation $(f1,4\bar{5})$ signifying the pair consisting of
the face $f1$ and the component given by the traversing the edge $e4$
followed by the reverse of the edge $e5$.

We will need the notion of a thickening of $G$ - by which we will
understand a sufficiently small neighbourhood of $G$ with respect to
some Riemannian metric on $S^2$. A moment's thought shows that there
is a natural bijection between the set of boundary components of such
a thickening of $G$ and what we earlier called $\widetilde{F}$.
A clockwise traversal of the boundary component corresponding to
$(f,c) \in \widetilde{F}$ (under the above  bijection) leads naturally
to what we would like to term a clockwise enumeration of $D_{(f,c)}$. 
Denote this enumeration by 
$\rho_{(f,c)} : \{1,\cdots,d_{(f,c)}\} \rightarrow D_{(f,c)}$.

In our example, the sets $D_{(f,c)}$, with their members listed in a
choice of such a clockwise order, are as follows:
\begin{eqnarray*}
D_{\tilde{f}1} &=& \{(e4,r),(e5,r)\}\\
D_{\tilde{f}2} &=& \{(e5,l),(e6,r),(e7,l)\}\\
D_{\tilde{f}3} &=& \emptyset\\ 
D_{\tilde{f}4} &=&
\{(e1,r),(e3,r),(e3,l),(e2,r), (e2,l),(e1,l)\}\\ 
D_{\tilde{f}5} &=& \{(e4,l), (e7,r),(e6,l)\}~.
\end{eqnarray*}

Now, $D(F)$ is the disjoint union of the $D_{(f,c)}$ as $(f,c)$ 
range over $\widetilde{F}$ and so the element $\otimes_{(f,c) \in 
\widetilde{F}}
\widetilde{\rho_{(f,c)}}(\delta^{-1}\Delta_{d_{(f,c)}}(h))$ is a well-defined element of
$H^{\otimes D(F)}$ which is independent of the choice of clockwise enumerations
of the $D_{(f,c)}$'s.

Finally, consider the map $\mu \circ (id \otimes S): H^{\otimes \{l,r\}} = 
H^{\otimes 2} \rightarrow H$.
In this, $\{l,r\}$ is mapped to $\{1,2\}$ by $l \mapsto 1$ and $r \mapsto 2$.
The tensor product map $\otimes_{e \in E} (\mu \circ (id \otimes S)) :
H^{\otimes D(F)} = H^{E \times \{l,r\}} \rightarrow H^E.$ 
Define $F(G,H)$ to be the image under this map of $\otimes_{(f,c) \in 
\widetilde{F}}
\widetilde{\rho_{(f,c)}}(\delta^{-1}\Delta_{d_{(f,c)}}(h))$.
The element of interest is $F(G,H^*)$ which is obtained by replacing $h$
by $\phi$ in the above expression. Explicitly, we have
\be \label{fgh*def}
F(G,H^*) = \delta^{\sigma(G)}\otimes_{e \in E} \phi^{L(e)}_{i(e)}
S\phi^{R(e)}_{j(e)} 
\ee
where  (i) $\sigma(G) = -|\widetilde{F}| + 2|\{v \in V: d_v = 0\}|$; and
(ii) $L,R,i,j$ are functions defined on $E$ and with appropriate
ranges, so that $(e,l)$ is the $i(e)$-th
element of $D_{L(e)}$ while $(e, 
r)$ the $j(e)$-th element of  of $D_{R(e)}$, for any edge $e \in E$.
(Thus, for example, $L,R :E \rightarrow \widetilde{F}$.)

In our example, $F(G,H^*) \in
(H^*)^{\otimes 7}$; and the prescription unravels to yield
\be \label{fgh*}
F(G,H^*) = \delta^{-3}\left(\phi^4_6S\phi^4_1 \otimes \phi^4_5
  S\phi^4_4 \otimes 
\phi^4_3S\phi^4_2 \otimes \phi^5_1S\phi^1_1 \otimes \phi^2_1S\phi^1_2
\otimes \phi^5_3S\phi^2_2 \otimes \phi^2_3S\phi^5_2\right)~.\ee

The remainder of this section is devoted  to proving the following:

\begin{Proposition} \label{grpr}
$$
F(G,H^*) =   F^{\otimes E}~(V(G,H))
$$
 for any spherical graph $G$.
\end{Proposition}

Our proof goes through the machinery of planar algebras but it
would be desirable to find a direct proof.

We use $G$ to construct a network in the Jones sense on $S^2$.
This network will be denoted $N = N(G)$. To construct $N$, choose
a thickening of $G$, as described above.
Colour this subset of $S^2$ black. Each edge of $G$ now appears
as a thin black band in this subset.
Replace this portion of the band by introducing a 2-box as indicated
 below:
\[\hpic {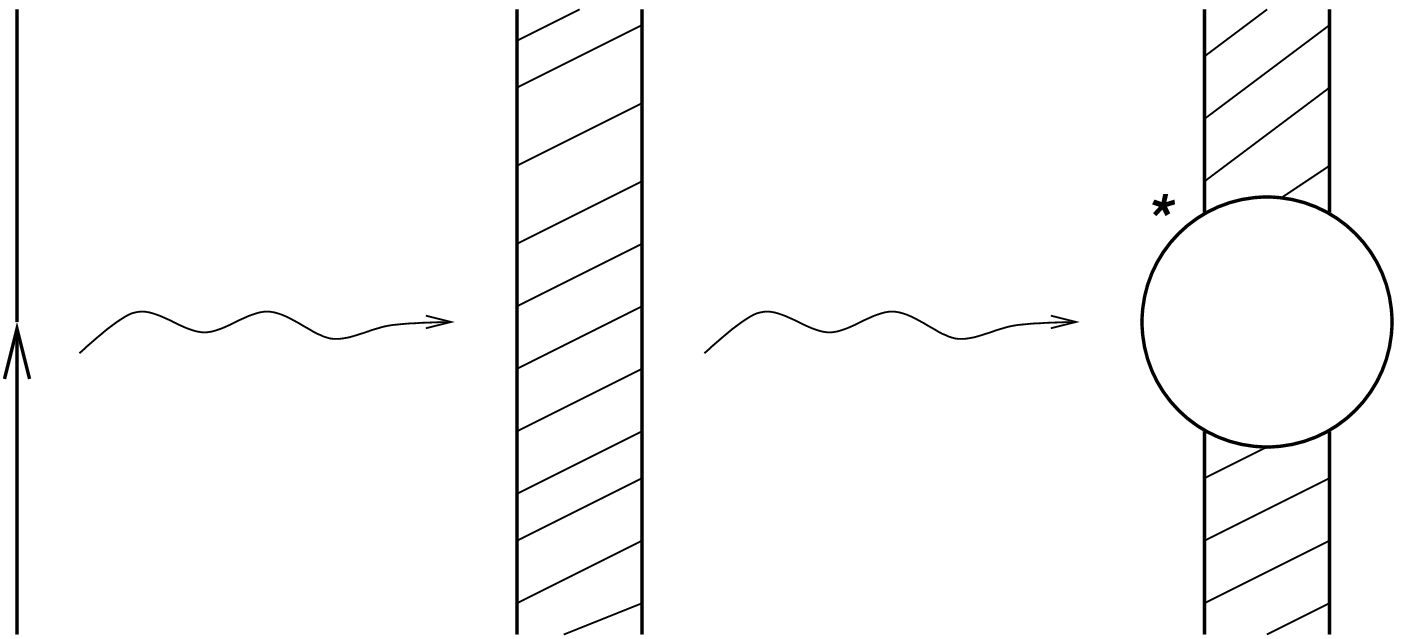} {1in} \]
with
the orientation of the edge determining the position of the $*$.
This yields our network $N$ on the sphere; note that $N$ has
only 2-boxes.
From the construction it should be clear that there are natural
bijections between the sets of black regions, white regions and 2-boxes
of $N$ and the sets of vertices, faces and edges of $G$ respectively. 

If $P$ is any spherical planar algebra, the partition function of $N(G)$
specifies a function from $(P_2)^{\otimes E}$ to $P_{0_+}$.
In particular, if $P = P(H)$,  this partition
function may be identified with a linear map from $H^{\otimes E}$ to
$k$ or equivalently, with an element of $(H^*)^{\otimes E}$.
We assert that this element is exactly $F(G,H^*)$. Explicitly, we need to verify that
\be \label{tpt}Z_{N(G)}(\otimes_{e \in E} a^e) = (F(G,H^*))(\otimes_{e
  \in E} a^e) ~\forall a^e \in H~.\ee
By definition of $F(G,H^*)$, we have
\begin{eqnarray*}
(F(G,H^*))(\otimes_{e \in E} a^e) &=& \delta^{\sigma(G)}\prod_{e \in
  E} \left( \phi^{L(e)}_{i(e)} S\phi^{R(e)}_{j(e)}\right)(a^e)\\
&=& \delta^{\sigma(G)}\prod_{e \in
  E}  \phi^{L(e)}_{i(e)}(a^e_1) \phi^{R(e)}_{j(e)} (Sa^e_2)\\
&=& \delta^{\sigma(G)}\prod_{Q \in \widetilde{F}} \left[  \left( \prod_{e \in
  E : L(e) = Q}  \phi^{Q}_{i(e)}(a^e_1)  \right) \left(   \prod_{e \in
  E : R(e) = Q}  \phi^{Q}_{j(e)} (Sa^e_2) \right) \right]\\
&=& \delta^{\sigma(G)}\prod_{Q \in \widetilde{F}} \phi^Q \left(
  \prod_{i=1}^{d_Q} T^Q_i a^{\rho_Q(i)}_{\e^Q_i} \right)\\
\end{eqnarray*}
where  $(T^Q_i, \e_Q(i)) = \left\{ \begin{array}{ll} (id,1) & \mbox{if
    } (\rho_Q(i),l) \in D_Q\\ (S, 2)  & \mbox{if } (\rho_Q(i),r) \in
    D_Q \end{array} \right. $~.

The proof of the asserted equation (\ref{tpt}) follows immediately
from Corollary 3 of \cite{KdySnd2}. (One only needs to note that the
`loops' of that 
prescription are in bijection with members of $\widetilde{F}$, and
exercise a little caution - in case $G$ has isolated vertices, so that
$N(G)$ has isolated loops - to see that the powers of $\delta$ also match.)

We next assert
that with identifications as above, $Z_{N^-} = V(G,H^*)$. This
assertion is proved exactly like the equation $Z_N = F(G,H^*)$ was proved
- after having observed that the black and white regions for the
network $N^-$, correspond to the white and black regions for $N$.

Applying Proposition \ref{main} to $N = N(G)$, 
\[F(G,H^*) = Z^{P(H)}_N = Z^{P(H^*)}_{N^-} \circ
F^{\otimes E} = V(G,H) \circ
F^{\otimes E} = F^{\otimes E}(V(G,H))~.\]
The first $V(G,H)$ is regarded as an element of $(H^{**})^{\otimes E}$ while 
the second is
regarded as an element of $H^{\otimes E}$, and the last equality follows from 
$x\circ F(y) = (F(x))(y)$.

So, Proposition \ref{grpr} has been finally proved.

We finally wish to observe a consequence of this proposition.

\begin{Corollary}\label{ngon}
In any semisimple cosemisimple Hopf algebra, we have 
\begin{eqnarray*}
(a) ~h^0_1 Sh^1_2 \otimes h^1_1 Sh^2_2 \otimes \cdots \otimes h^{(n-1)}_1
Sh^0_2 &=& \delta^n F^{\otimes n} (\Delta_n \phi)\\
(b) ~h^1_1 Sh^0_2 \otimes h^2_1 Sh^1_2 \otimes \cdots \otimes h^0_1
Sh^{(n-1)}_2 &=& \delta^n F^{\otimes n} (\Delta^{op}_n \phi)
\end{eqnarray*}
for any $n \geq 1$.
\end{Corollary}

To prove this, consider the special case of Proposition
\ref{grpr} corresponding to $G$ being a cyclically oriented
$n$-gon. Write $V = \{0,1,\cdots ,n-1\}, E = \{e0,e1,\cdots ,e(n-1)\}, F =
\{in, out\}$, and make `cyclically symmetric' choices as below (where
we illustrate the case $n=6$:
\[\hpic {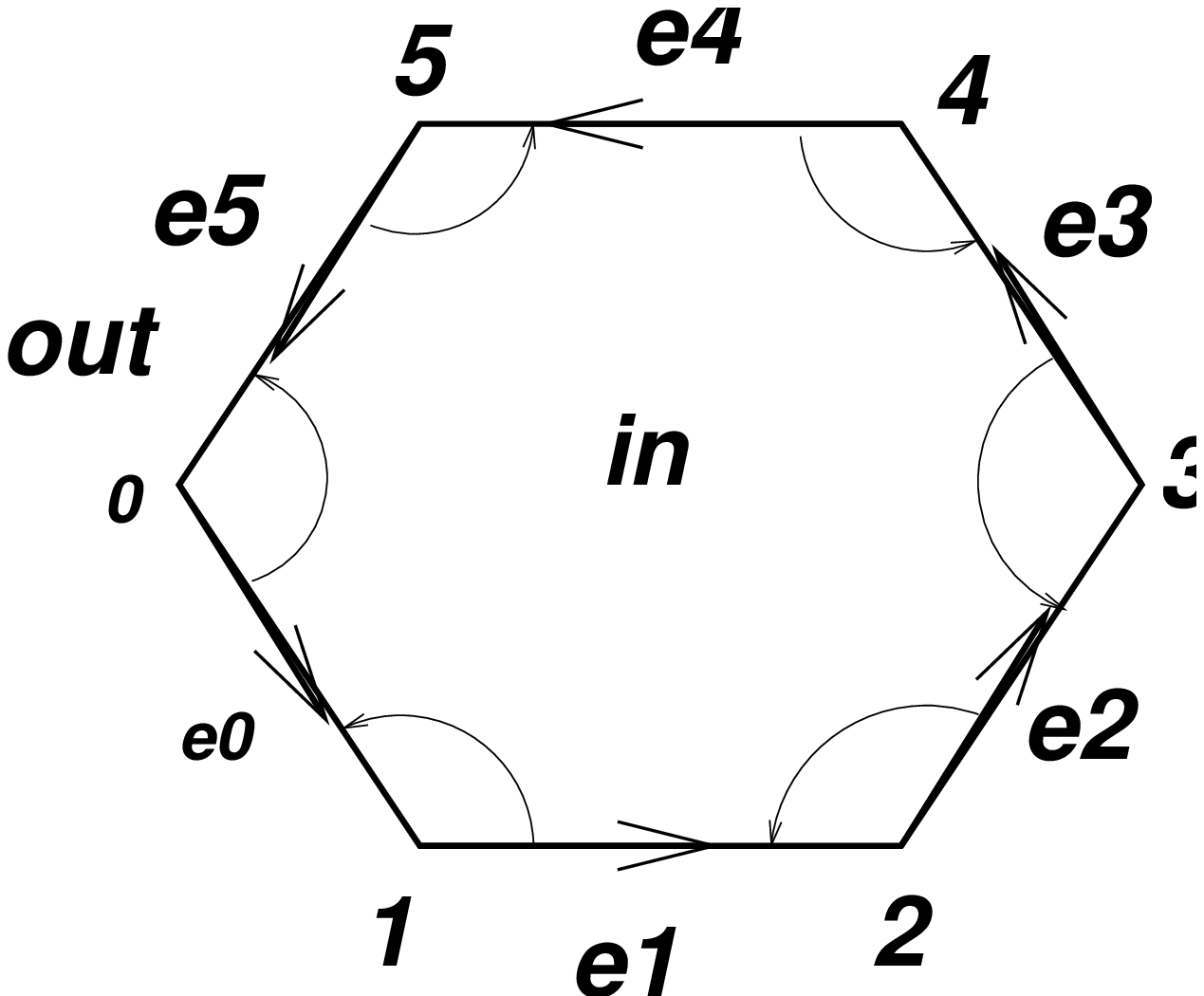} {2in} \]
We set
\[ D_i=\{(ei,0),(e(i-1),1)\}, ~\forall ~0 \leq i < n ~,\]
with addition modulo $n$.
Further, $\widetilde{F} = F$, and we choose
\[D_{in} = \{(e(n-1),l), \cdots ,(e1,l), (e0,l)\} \mbox{ and } 
D_{out} = \{(e0,r), \cdots , (e(n-1),r)\}~.\]
Our prescriptions yield
\[V(G,H) = \delta^{-n}(h^0_1 Sh^1_2 \otimes h^1_1 Sh^2_2 \otimes
\cdots \otimes h^{(n-1)}_1Sh^0_2) ~;\]
and
\[F(G,H^*) = \delta^{-2}(\phi^{in}_n S\phi^{out}_1 \otimes \phi^{in}_{n-1}
S\phi^{out}_2 \otimes \cdots \otimes \phi^{in}_1S\phi^{out}_n) ~.\]

Since $S^{\otimes n}(\Delta_n(a)) = \Delta_n^{op}(Sa)$ in any Hopf
algebra, this simpliies to 
\begin{eqnarray*}
F(G,H^*) &=& \delta^{-2} \Delta_n^{op}(\phi^{in}S\phi^{out})\\
&=& \Delta_n^{op} (\phi)~,
\end{eqnarray*}
the final equality being a consequence of the fact that
$\phi^2 = \delta^2 \phi$ and $S\phi = \phi$.

So we deduce from Proposition \ref{grpr} that
\[F^{\otimes n}\left(\delta^{-n}(h^0_1 Sh^1_2 \otimes h^1_1 Sh^2_2 \otimes
\cdots \otimes h^{(n-1)}_1Sh^0_2)\right) = \Delta_n^{op} (\phi)~;\]
and since $F^{-1} = F\circ S$, we conclude that
\begin{eqnarray*}
h^0_1 Sh^1_2 \otimes h^1_1 Sh^2_2 \otimes \cdots \otimes h^{(n-1)}_1Sh^0_2
&=& \delta^{n} (F\circ S)^{\otimes n} (\Delta_n^{op}
(\phi))\\
&=& \delta^{n} F^{\otimes n} (\Delta_n(\phi))~,
\end{eqnarray*}
thus establishing (a).
By applying $S^{\otimes n}$ to both sides of (a), (b) follows.


\begin{thebibliography}{amsalpha}

\bibitem[BhmNllSzl]{BhmNllSzl} Gabriella Bohm, Florian Nill and Kornel Szlachanyi,
{\em Weak Hopf algebras I. Integral theory and $C^*$-structure}, Journal of Algebra, 221, (1999) 385-438.

\bibitem[BrrWst]{BrrWst} Barrett, John W.; Westbury, Bruce W. {\em The
  equality of $3$-manifold invariants},  Math. Proc. Cambridge
  Philos. Soc.  118  (1995),  no. 3, 503--510.  

\bibitem[DttKdySnd]{DttKdySnd1} Sumanth Datt, Vijay Kodiyalam and V. S. Sunder,
{\em Complete invariants for complex semisimple Hopf algebras},
Math. Res. Lett., 10, (2003) 571-586.

\bibitem[Jns]{Jns} V. F. R. Jones, {\em Planar algebras I}, New
  Zealand J. of Math., to appear. e-print arXiv : math.QA/9909027

\bibitem[KdyLndSnd]{KdyLndSnd} Vijay Kodiyalam, Zeph Landau and V. S. Sunder,
{\em The planar algebra associated to a Kac algebra}, Proc. Indian Acad. of
Sciences, 113, (2003) 15-51.

\bibitem[KdyPtiSnd]{KdyPtiSnd} Vijay Kodiyalam, Vishwambhar Pati and V. S.
Sunder, {\em Subfactors and 1+1 dimensional TQFTs}, Preprint.

\bibitem[KdySnd1]{KdySnd1} Vijay Kodiyalam and V. S.
Sunder, {\em On Jones' planar algebras}, J. Knot theory and its ramifications,
13, (2004) 219-247.

\bibitem[KdySnd2]{KdySnd2} Vijay Kodiyalam and V. S. Sunder, {\em
The planar algebra of a semisimple cosemisimple Hopf algebra},
e-print arXiv math.QA/0506153.

\bibitem[Kpr]{Kpr} G. Kuperberg, {\em Involutory Hopf algebras and
3-manifold invariants},
Internat. J. Math., 2, (1991) 41-66.
 
\bibitem[Lnd]{Lnd} Zeph Landau, {\em Exchange Relation Planar
    Algebras}, Proceedings of the Conference on Geometric and
    Combinatorial Group Theory, Part II (Haifa, 2000).  Geom. Dedicata
    95  (2002), 183--214.  

\bibitem[PrsSss]{PrsSss} V.V. Prasolov and A B. Sossinsky, {\em Knots,
    Links, Braids and 3-Manifolds : Introduction to  the New
    Invariants in Low-Dimensional Topology : (Transl. of
    Math. Monographs 154)}, Amer. Math. Soc., USA, 1997.

\end{thebibliography}
\end{document}